\documentclass[11pt]{article}
\usepackage{amsfonts, amssymb, amsmath, amscd, longtable, tabularx, multicol}
\usepackage{graphicx}
\newtheorem{theorem-intro}{Theorem}[section]
\newtheorem{theorem}{Theorem}
\newtheorem{lemma}{Lemma}

\newtheorem{corollary-intro}[theorem-intro]{Corollary}
\newtheorem{remark}{Remark}

\newtheorem{example}{Example}

\def\be{\begin{enumerate}}
\def\ee{\end{enumerate}}
\def\bi{\begin{itemize}}
\def\ei{\end{itemize}}
\def\xi{\mathcal X}

\title{PA mapping classes with minimum dilatation and Lanneau-Thiffeault polynomials }
\author{Joan S. Birman}
\date{}

\begin{document}
\maketitle

\begin{abstract}
It has been known since 1981 that if one fixes an orientable surface $S$ of genus $g$, then there is a real number $\lambda_{min,g} > 1$ that is the dilatation of a pA diffeomorphism of $S$, and every other pA diffeomorphism of $S$
 has  dilatation $\geq \lambda_{min,g}$.  We will show how a little-known theorem about digraphs gives some insight into $\lambda_{min,g}$.
\end{abstract}

\section{Introduction}  \label{S:introduction}
Let $S = S_g$ be a closed, oriented surface of genus $g\geq 2$.  The {\em mapping class group} ${\rm Mod}(S)$ is the  group
$\pi_0({\rm Diff^+}{S})$, where admissible diffeomorphisms preserve orientation.  
In \cite{Th} Thurston proved a far-reaching classification theorem for elements in Mod($S)$.  He showed that for each homotopy class in Mod$(S)$, there is a diffeomorphism $f$   such that $f$ is one of three mutually exclusive types:  {\it finite order}, {\it reducible} or (the generic case \cite{R}) 
 {\em pseudo-Anosov}.   In the pseudo-Anosov or pA case, even though $f$ and all its powers fix no simple closed curve on $S$,  Thurston discovered rich structure: there exists
a pair of $f$-invariant transverse measured
foliations ${\mathcal F}^u, {\mathcal F}^s$  of $S$ and a real number $\lambda(f)$, the
{\em dilatation} of $f$, such that $f$ maps leaves of ${\mathcal F}^u, {\mathcal F}^s$ to leaves and multiplies the measure on ${\mathcal
F}^u$ (resp. ${\mathcal F}^s$) by $\lambda$ (resp.  $\frac1\lambda$).   

The real number $\lambda$ is positive, it's greater than 1, and it's a conjugacy class invariant of $f$ in Mod$(S)$.  It's  the largest real root of the characteristic polynomial $\xi_T(x)$ = det$(xI - T)$ of a non-negative integer matrix  $T= (t_{i,j})$.  The matrix  $T$ arises via the action $f_\star$ of $f$ on a graph $G$ that was introduced and studied by Bestvina and Handel in \cite{BH}.  Let $\{e_1,\dots,e_m\}$ be the edges of $G$. Then $f_\ast$ acts on the edges, sending vertices to vertices and sending each edge $e_i$ to an edge sequence $f_\star(e_i) = e_{i_1}\dots\_{i_{k_i}}$.  The entry $t_{i,j}$ in $T$ is the number of times that $f_\star(e_i)$ crosses $e_j$, in either direction.  The matrix $T$ is {\it primitive}, that is there is a positive integer $k$ such that every entry in $T^k$ is strictly positive \cite{FLP}, which implies that $T$ is irreducible.   
The Perron-Frobenius theorem \cite{Se} applied to $T$ then shows that  $\xi_T(x)$ has a unique largest real root.  This root is $\lambda(f)$.     It was proved by Arnoux-Yoccoz \cite{AY} and Ivanov  \cite{Iv}  that for each fixed choice of surface $S$, there is a lower bound $1< \lambda_{{\rm min}, g}$ for the dilatation of all pseudo-Anosov (pA) maps on $S$, and moreover this lower bound is achieved for some pA mapping class $f$.  Our work in this paper concerns the somewhat mysterious number  $\lambda_{{\rm min}, g}$.   

Interest in $\lambda_{min,g}$ has implications that go beyond surface mappings to 3-manifolds. Going back to the original work of Thurston, recall that he proved that the mapping torus of a surface diffeomorphism $f$ is hyperbolic if and only if $f$ is pA.  This suggests that questions about  pA maps with very low dilatation and hyperbolic fibered $M^3$ with very low volume might be 
related.  Computations, and recent work on hyperbolic volume (e.g. see \cite{FLM}) shows that when the dilatation is small, a relationship exists.  Examples given in \cite{AD} make the relationship concrete in special cases.  Thus it's an interesting open problem to find $\lambda_{min,g}$.   

Upper bounds for $\lambda_{min,g}$  have improved over the years as new examples of surface diffeomorphisms with low dilatation were discovered, and of course each new example with lower dilatation than previously known examples becomes a new candidate for $\lambda_{min,g}$.   
While the known upper bounds (see \cite{AD}) were due to several different authors taking different approaches to the problems, a striking fact stood out and became the takeoff point for the work in this paper:
The family of integer polynomials:    
\begin{equation} \label{E:LT}
\xi_{d,a}(x) = x^{2d} - x^{2d-a} - x^d - x^a + 1, \ \ 1\leq a\leq d-1,  
\end{equation}
are known as {\it Lanneau-Thiffeault polynomials} because they were first singled out  as a class in \cite{LT}.  Let $\lambda_{d,a}$ be the unique largest real root $>1$ of $\xi_{d,a}(x)$.  In \cite{H} Hironaka proved that if $g\geq 5$ then:
\begin{eqnarray} 
 \label{E:H1} \nonumber
\lambda_{min,g} &\leq & \lambda_{g+1,g} \ \ \ \ \ {\rm if} \   g = 7+3n \ {\rm or} \  8+3n, \ n\geq 0, \\ 
\label{E:H}
& \leq & \lambda_{g+1,g-2}\ \  \ \ \    {\rm if} \ g=6+3n, \ n\geq 0. 
\end{eqnarray}

Let $f_{d,a}$ be a pA map whose dilatation is $\lambda_{d,a}$.  By \cite{H} such an $f_{d,a}$ exists.  We noticed that the transition matrix $T_{d,a}$ is always a (0,1)-matrix having precisely $2d+2$ non-zero entries.    This seemed interesting to us.  If a pseudo-Anosov map has very low dilatation, and if $\{e_1,\dots,e_m\}$ are the edges of the Bestvina-Handel graph, then the edge paths $f_\star(e_i)$ must be short, because the growth rate is small. This is consistent with $T$ having a small number of non-zero entries.   Since $T$ has dimension $m$, it must have at least $m$ non-zero entries, so it is natural to define the {\it complexity} $c = c(T)$ of  $T$ to be the sum of the entries in $T$ minus $m$.    

There were also strong hints in the existing data  that a unifying concept that might cast light on the pA maps of least dilatation  could be the fact that $T$ is non-negative.   With this last fact in mind, we were lead to think about an alternative method of describing $T$.   A {\it digraph}  is a directed graph.  Digraphs enter into our work because one may associate to every transition matrix its {\it associated digraph}. 
If $T$ has dimension $m\times m$, then the digraph $D = D_T$ associated to $T$ has $m$ vertices, and for every pair $(i,j), \ 1\leq i,j \leq m$ it has $t_{i,j}$ directed edges from vertex $v_i$ to vertex $v_j$.  Clearly $T$ and $D$ determine one-another.  
The matrix $T$ is irreducible if and only if  $D$ is {\it strongly connected}, that is there is a directed path from every vertex $v_i$ of $D$ to any other vertex $v_j$.   The matrix $T$ is primitive if and only if $D$ is strongly connected and the gcd of all the path lengths is 1. Define $\xi_{D_T}(x) = \xi_T(x)$. 
We will introduce a tool that, to the best of our knowledge, has not been applied before this to study the structure of $ \xi_T(x)$ when $T$ is topologically induced by a pA map.  The {\it Coefficient Theorem for Digraphs}  (see Theorem~\ref{T:coef} below) tells us that  there is structure in $D$ that is not easily detected in $T$, and this proved to be very very useful.  

Our main result is: 

\begin{theorem} \label{T:LT}  
Let $T$ be the transition matrix for the action of a pA map $f$ on a surface $S$ of genus $g$ that realizes $\lambda_{min,g}$.  Let $\xi_T(x)$ be its characteristic polynomial and let $c$ be its complexity. 
\begin{enumerate}
\item [\rm{(i)}] If $c\leq 2$, then  the unique possibility is that $c=2$ and $\xi_T (x)$ is the LT-polynomial $\xi_{d,a}$ defined in {\rm (\ref{E:LT})} above, for some  $d$ with $g\leq d\leq 3g-3$.   
\item [\rm{(ii)}]   If $c\leq 5$, then there is precisely one additional possibility:  $\xi_T(x)$ belongs to the following family of examples of complexity $4$ : 
\begin{equation} \label{E:c=4 examples}
 \xi_{d,\vec{a}}(x) = x^{2d} -\sum_{i=1}^4 (x^{2d-a_i} +x^{a_i})+ \sum_{j=1}^3 (x^{2d -a_1-a_j} +x^{a_1+a_j}) -x^d + 1,
 \end{equation}
where $\vec{a}=(a_1,a_2,a_3,a_4)$, each $a_i$ an integer  $\geq 2$ and $\sum_{i=1}^4a_i = 2d$. 
\end{enumerate}
\end{theorem}
\begin{remark}
{\rm When the work in this paper began we had been studying, together with Brinkmann and Kawamuro, the factorization properties of $\xi_T(x)$ for arbitrary pA mapping classes.  It is proved in \cite{BBK}  that, with our assumptions, 
$\xi_T(x)$ can be assumed to be either palindromic or anti-palindromic.  That is, 
\begin{eqnarray}\label{E:coef}
\xi_T(x) &=& x^m + \sum_{i=1}^m b_ix^{m-i}, \ {\rm where}  \\
\nonumber
b_m &=& +1 \ {\rm and} \  b_i = b_{m-i} \  {\rm if} \ 1 \ \leq i \leq m-1  \ {\rm or} \\
\nonumber
b_m &=& -1 \ {\rm and} \  b_i = - b_{m-i} \  {\rm if} \ 1 \  \leq i \leq m-1.
\end{eqnarray}
This gave us structure that had not been used before as a tool in the study of the dilatations of low dilatation pA maps.   
As will be seen, the Perron-Frobenius matrices that can be realized topologically, as transition matrices for a pA map,  are a very small subset of all Perron-Frobenius matrices of the same dimension. }
\end{remark}


{\bf Acknowledgements}  We gratefully acknowledge stimulating and helpful conversations with  Francesco Belardo, Richard Brualdi, Dan Margalit and Andrew Putman, Eriko Hironaka,  Michael Polyak,  Igor Rivin, Saul Schleimer,  Adam Sikora.   

\section{Proof of Theorem~\ref{T:LT}}\label{S:the proof}  

Essentially all of our work in the proof of Theorem~\ref{T:LT} will be done in the setting of the digraph $D=D_T$, and so we begin with the statement of a theorem about digraphs which will play a crucial role in the proof. 

A subdigraph of $D$ is a {\it cycle} if every vertex has in-valence and out-valence 1.  A subdigraph of $D$ is {\it linear} if it is a union of pairwise disjoint cycles.   The  symbol $\mathcal L_i$ will be used to denote the set of all $i$-vertex linear subdigraphs $L_{i,j}$ of a $D$.  If there is only one,  then we write $L_i$ instead of $L_{i,j}$.  The symbol $n(L_{i,j})$ denotes the number of cycles in $L_i$.  The following theorem tells us that the characteristic polynomial of a digraph (and so also of the associated transition matrix) is determined by its linear subdigraphs:
\begin{theorem}\label{T:coef} {\bf (The Coefficient Theorem for Digraphs, or the CT).} See {\rm \cite{CDS}.} 
Let $D$ be a digraph with $m$ vertices and let $\mathcal L_i = \cup _jL_{i,j}$ be the set of all linear subdigraphs of $D$ having precisely $i$ vertices.  For each $L_{i,j}\in \mathcal L_m$, let $n(L_{i,j})$ denote the number of cycles in $L_{i,j}$.  Then
the characteristic polynomial of $D$ is
\begin{eqnarray} \label{E:characteristic poly}
\xi_D(x) & = & x^{m} +\sum_{i=1} ^{i=m} b_ix^{m-i}, \ \ {\rm where} \\
\nonumber
 b_i & = &  
\sum_{L_{i,j} \in \mathcal L_i} (-1)^{n(L_{i,j})}
\end{eqnarray}
\end{theorem}

\begin{example}  \label{Ex:characteristic polynomial of D(2,2), (d,a)=(7,6)} 
{\rm We illustrate the CT with an example.  The digraph $D$ in  sketch (a) of 
Figure 1  
has 14 vertices, labeled $1,2,\dots,14$, and 16 edges. 
It is strongly connected.   Vertices $1$ and $14$ (resp. $1$ and $9$) have out-degree (resp. in-degree) 2, whereas all others have out-degree (resp. in-degree) 1. 
 \label{F:digraph1}
\begin{figure}[htpb!] 
\centerline{\includegraphics[scale=.50] {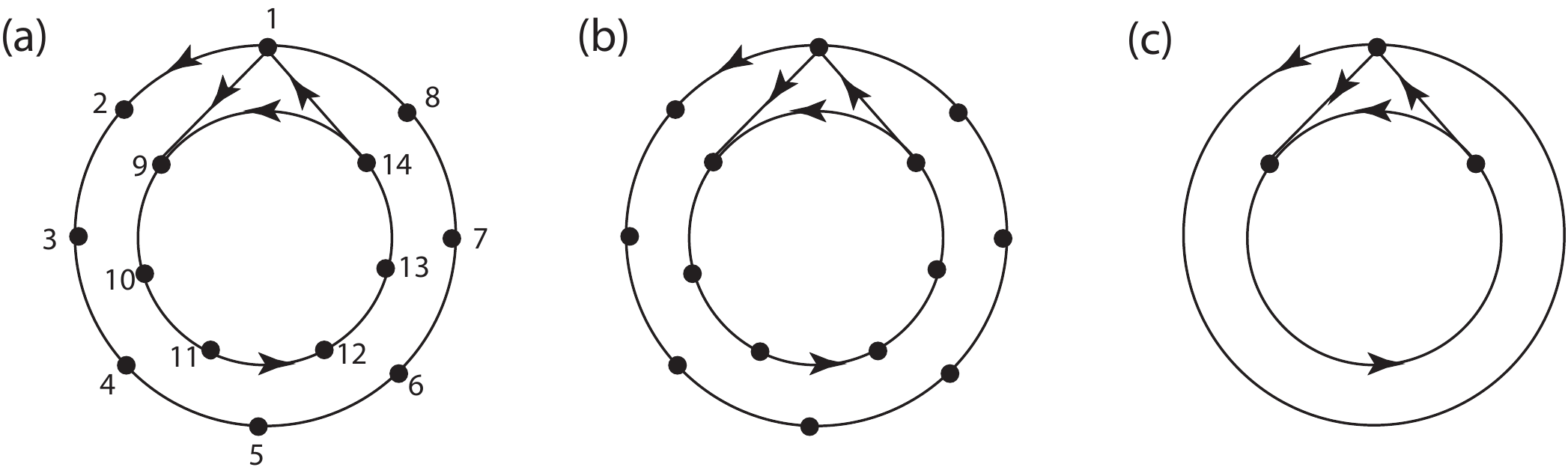}}
\caption {(a) A strongly connected digraph $D$ with 14 vertices and 16 edges; (b) the same digraph without vertex labels; (c) The {\it shape} of $D$ }
\end{figure}
Sketch (b) shows the same digraph without a choice of numbers on the vertices.  It determines all matrices that can be obtained from $T$ after an arbitrary permutation of the rows and columns.  Sketch (c) shows the {\it shape} of $D$, and so incorporates all the information for an entire class of digraphs with varying numbers of vertices.  In particular, the number of cycles, and the number of vertices that have in-degree or out-degree $>1$ can be seen from sketch (c).  The pictures in sketches (b) and (c), but not (a), will be used in this paper.  That is one reason why $D$ reveals structure that may be very hard to see in $T$.

We compute the characteristic polynomial of $D$, using sketch (b).   It has four linear subdigraphs, with 6,7,8 and 14 vertices. 
The first three contain one cycle and the fourth contains two disjoint cycles, therefore $n(L_6)=n(L_7)=n(L_8)=1$, whereas $n(L_{14})=2$. Therefore the coefficients $b_6 = b_7 = b_8 = -1$ and $b_{14} = +1$.   There are no additional linear subgraphs, so that all other $b_i$ are 0.    The Coefficient Theorem tells us, without any calculation, that the characteristic polynomial of this digraph and so also of the associated matrix, is the palindromic polynomial  $\xi_{7,6}(x) = x^{14} - x^8 - x^7 -x^6 +1$, an LT-polynomial.    Note that this is independent of the way that the vertex labels are assigned, since all we used was sketch (b).   Indeed, sketch (b) shows that $D$ determines not just $T$, but all of the matrices obtained from $T$ by an arbitrary permutation of the rows and columns of $T$.  This example will be discussed again in Example~\ref{Ex:LT} (see the bottom leftmost sketch in Figure 3).  We will prove, very soon, that there are 5 possible digraphs with 14 vertices that have this polynomial.  The polynomial occurs in \cite{AD} as the monodromy map of a fibered hyperbolic 3-manifold of very low volume.  The associated transition matrix could be any one of the 5 possible digraphs in the bottom row in Figure 3. It has the lowest known dilatation for all pA maps on a closed orientable surface of genus 5.

As for sketch (c), let $A_0, A_1,A_2$ be the 3 cycles in sketch (c), where $A_0$ (resp. $A_1$) is the outermost (resp. innermost) cycle in the sketch,  and $A_2$ is the cycle that includes the connecting edges.  Thus $A_0\cap A_1 = \emptyset$ but $A_0\cap A_2 = 2$ vertices, $A_1\cap A_2$ = 1 vertex. Let $a_0, a_1,a_2$ be the number of vertices in $A_0, A_1,A_2$.   Sketch (c) describes an entire class of polynomials, namely the polynomials 
$x^m -x^{m-a_0} - x^{m-a_1} -x^{m-a_2} + x^{m-a_0-a_1}$.   Thus there is a great deal to be learned by looking at the shape of a digraph, without specifying all of its vertices.   \ \ \ \ \ \ \ \ \ \ \ \ \ \ \ \   $\square$

}\end{example}



Our plan for the proof of Theorem~\ref{T:LT} is to study the possibilities for $\xi_D(x)$ under the given hypotheses, which as we shall see are very restrictive. 

\begin{lemma} \label{L:BBK}
Let $\xi_D(x)$ be the characteristic polynomial of $D$.  Then:   
\begin{enumerate}
\item [{\rm (1)}]  $\xi_D(x)$ is palindromic or anti-palindromic. In particular,  $b_m = \pm 1$.  
\item [{\rm (2)}]  $D$ has a linear subdigraph $L_m$ which uses the $m$ vertices in $D$ and is a disjoint union of $n$ cycles. 
\item [{\rm (3)}]  $n\leq c$.
\item[{\rm (4)}]  If $\xi_D(x)$ is antipalindromic, then $\xi_D(1) = 0$.
\end{enumerate} 
\end{lemma}
{\bf Proof:}  

(1) The fact that $\xi_D(x)$ is palindromic or anti-palindromic  is proved  in \cite{BBK}.   

 (2) From (1) it follows that $b_m =  \pm 1$, where $m$ is the dimension of $T$ or, equivalently, the number of vertices in $D$.  By the  CT it then follows that $D$ contains at least one  linear subdigraph $L_m$ which uses all $m$ vertices of $D$.   This subdigraph, being linear,  is a union of, say, $n \geq 1 $ disjoint cycles. 
 
(3) The complexity $c$ of $D$ is the number of edges minus the number of vertices.  The subdigraph $L_m$ has $m$ vertices, and since a cycle has the same number of vertices as edges, it also has $m$ edges.  But $D$ is strongly connected,  so it must have at least $n$ connecting edges, so the total number of edges in $D$ is $\leq n+m$.   The number of vertices in $D$ is $m$.  Therefore the complexity $c\geq n$.     

(4) The terms in an antipalindromic polynomial must occur in pairs with opposite coefficients, so if $\xi_D(x)$ is antipalindromic, then  the sum of the coefficients, that is $\xi_D(1)$, must be 0.  \ \ \ \ \ \ \ \ \ \ \ \ \ \ \ \ {$\square$}

We begin the proof of  Theorem~\ref{T:LT}.   By (2) of Lemma~\ref{L:BBK}, the digraph  $D(n,c)$  contains at least one linear subdigraph $L_m$ with $m$ vertices and $n$ cycles.    By (3) of Lemma~\ref{L:BBK}, we know that  $1\leq n \leq c$.   We may then parametrize the digraphs to be studied by the integer pairs $(n,c)$, ordered lexicographically.  Part (i) of Theorem~\ref{T:LT} relates to the cases when $(n,c) = (1,1),(1,2)$ or $(2,2).$   

If $(n,c) = (1,1)$,  then the linear subdigraph $L_m$ of part (2) of Lemma~\ref{L:BBK} has one oriented cycle, say $A$, with $m$ vertices.  See Figure 2.  There is also one additional oriented edge $e$ between two not-necessarily-distinct vertices of $A$.  Irrespective of the orientation of $e$, it follows that $D$ has a second cycle $A_1$ of length, say, $a_1$, where $1\leq a_1\leq m$.  There are no other linear subdigraphs. By the CT, it follows that $\xi_D(x) = x^m -x^{m-a_1} -x^{m-m}$, which implies, irrespective of the choice of $m$ or $a_1$,  that $\xi_D(1) = -1$ for any choice of $a_1$, contradicting (4) of Lemma~\ref{L:BBK}.
\begin{figure}[htpb!] 
\centerline{\includegraphics[scale=.5] {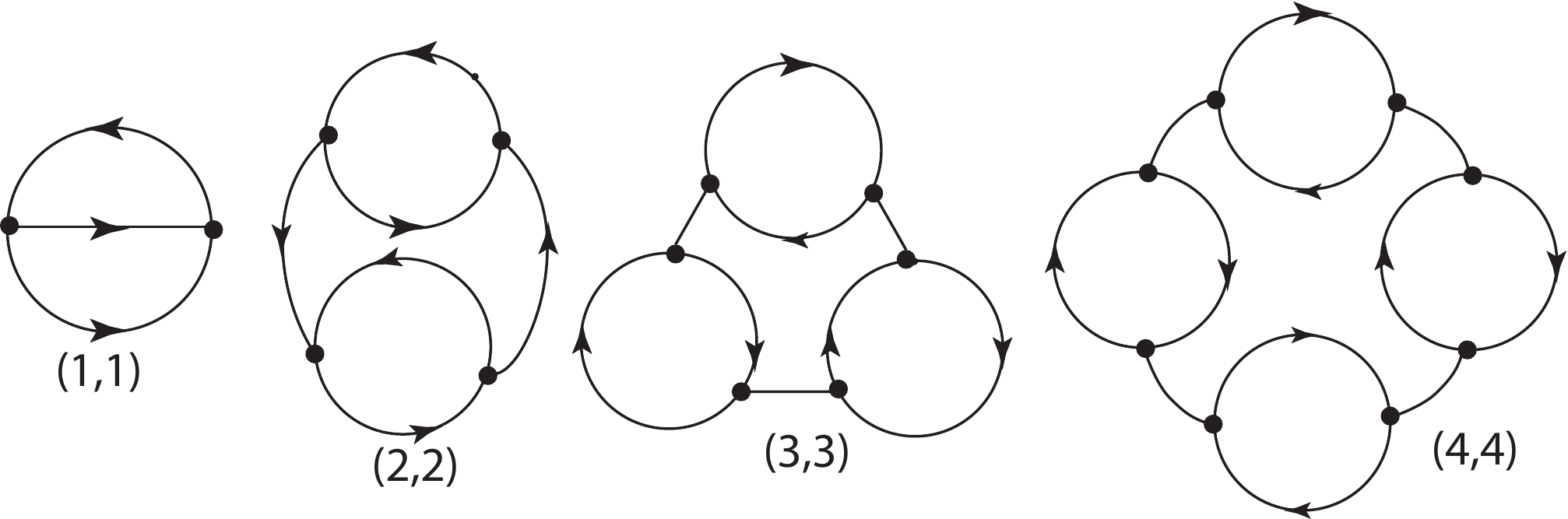}}
\caption{The shape of $D$ when $n=c =1,2,3,4.$ }
 \label{F:digraph2}
\end{figure}

If $(n,c) = (1,2)$ then the digraph we considered in the case $(n,c)=(1,1)$ is modified by the addition of a second oriented edge.  The two edges determine two new cycles   $A_1,A_2$, with $a_1,a_2$ vertices, also both $A_1$ and $A_2$ intersect $A$ non-trivially. There are two cases: 

Case (1):  $A_1\cap A_2 \not= \emptyset$.   
The three linear subdigraphs in $D$ each contain exactly one cycle,  that is $A, A_1, A_2$. 
$\xi_D(x) = x^m -x^{m-a_1} -x^{m-a_2} -x^{m-m} $, which implies that $\xi_D(1) = -2$, contradicting (4) of Lemma~\ref{L:BBK}.

Case (2): $A_1\cap A_2 = \emptyset$.  The only change from case (1) is that there is a fourth linear subdigraph in $D$, namely  $A_1\sqcup A_2$, and 
so $\xi_D(x) = x^m -x^{m-a_1} -x^{m-a_2} +x^{m-a_1-a_2} -1$.  But then $\xi_D(1) = -1$,  again contradicting (4) of Lemma~\ref{L:BBK}.

We turn to the case $(n,c) = (2,2)$, that is $L_m$ is a union of two disjoint oriented cycles, say $A_1$ and $A_2$, and $D$ consists of these cycles joined up by two edges, with $\vec{e_1}$ (respectively $\vec{e_2}$) directed from $A_1 \to A_2$ (resp. $A_2\to A_1$).  There are now 4 linear subdigraphs, made up from cycles $A_1, A_2, A_1\sqcup A_2$ and $A_3$ which includes the two connecting edges and paths on both $A_1$ and $A_2$.  See Figure 2 again.  This is true for every possible choice of orientations.  Observe that $a_1+a_2=m$.   By the CT, 
\begin{eqnarray}
\nonumber
\xi_D(x) &=& x^m -x^{m-a_1} -x^{m-a_2} -x^{m-a_3} + x^{m-a_1-a_2}\\
\label{E:parameters}
&=&  x^m -x^{m-a_1} -x^{a_1} -x^{m-a_3} + 1
 \end{eqnarray}
 It is immediate that $\xi_D(x)$ is palindromic if and only if $m=2d$ and $a_3 = m/2 = d$.  Simplifying notation by setting $a_1=a$,  one obtains the LT polynomial $\xi_{d,a}(x)$ that was defined in (\ref{E:LT}).   
  See Example \ref{Ex:LT} in the next section for a discussion of the finitely many possibilities. 
 
The proof of part (ii) of Theorem~\ref{T:LT} is similar.  Arguments just like the ones used to prove part (i) show that the only possibility is $(n,c) = (4,4)$,  yielding the generalized LT polynomial that is given in equation (\ref{E:c=4 examples}).  The digraph associated to (4,4) is illustrated in Figure 2.   In fact, the cases $(2k,2k), k\geq 1$ yield polynomials that are, in every case, easy generalizations of the LT polynomials. 

The proof that the cases $(n,c) = (2k+1,2k+1), k\geq 1$ cannot be realized topologically is a small generalization of  the proof for $(n,c)=(1,1)$.  

The cases $(n,c) = (2,k), k\geq 2$ are eliminated by using a fact that we have not used up to now:  The Perron-Frobenius Theorem, proved in \cite{Se}, asserts that if $D$ is primitive, then the Perron-Frobenius eigenvalue can never be reduced by adding more edges, therefore if there is an underlying linear subdigraph which looks like the ones in Figure 2 and uses all of the vertices, the spectral radius never decreases.  

The only case, up to $(n,c)\leq (5,5)$, where we had to consider more than the shape of $D$, and actually show that there was no way to place the vertices to obtain a palindromic polynomial, was $(n,c)=(4,5)$.   

There are cases where the combinatorics get out of hand.  If $(n,c) = (3,4)$, then  a few moments thought will convince the reader that the unique digraph with shape $(3,3)$ need not be a subdigraph using all vertices, and this situation persists as $n$ is increased.  Thus eliminating some of the cases required a lengthy calculation, with some computer help,  and without a really good upper bound on $c$ (a matter that is discussed in Example 6  below) the details are lengthy and unenlightening, therefore we omit them.  $\square$


\section{Examples, Conjctures, Open Problems} \label{S:Examples} 

\begin{example}  \label{Ex:LT}
{\rm We illustrate, via an example, what we already knew from the proof of Theorem~\ref{T:LT}: every LT polynomial occurs as the characteristic polynomial of a primitive Perron-Frobenius integer matrix, also the characteristic polynomial determines the digraph $D_T$ or, equivalently, the matrix $T$, up to a finite set of choices.  See Figure 3.  The example is the sequence of degree 14 polynomials $\xi_{7,a}, \ 1\leq a\leq 6.$ }
\begin{figure}[htpb!]  \label{F:digraph3}
\centerline{\includegraphics[scale=.55] {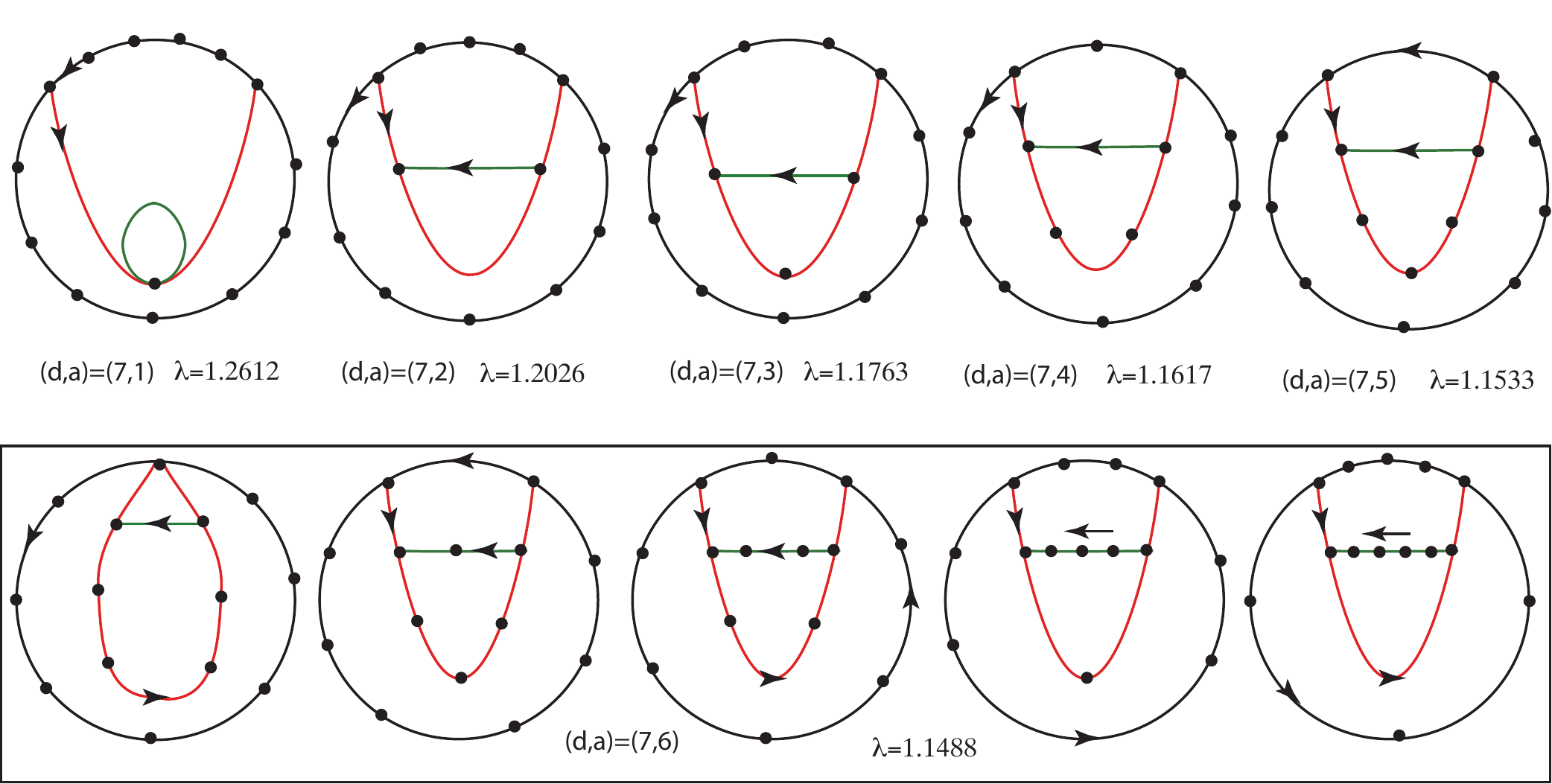}}
\caption{Digraphs whose characteristic polynomials are $\xi_{7,a}$, where $a=1,\dots,6$.  The boxed sketches are 5 distinct digraphs that realize $\xi_{7,6}$.}
\end{figure} 
\end{example}
{\rm  The boxed sketches show 5 distinct digraphs whose characteristic polynomials are all 
$$\xi_{7,6}(x) = x^{14}-x^8-x^7-x^6+1,$$ 
representing the 5 different choices of the parameter that we called $a_3$ in Equation (\ref{E:parameters}).  The leftmost one was discussed earlier in Example~\ref{Ex:characteristic polynomial of D(2,2), (d,a)=(7,6)}.  There are similar choices in the other 5 cases.

It's interesting  that $\xi_{7,6}(x)$ is the characteristic polynomial associated to the genus 5 example of Aaber-Dunfield in \cite{AD}.  Its dilatation is the minimum dilatation known, at this writing, that is realizable for genus 5.
In an e-mail interchange with Dunfield, he wrote: ``The structure of the invariant foliations is as follows.Ê There are 5 singularities with 3 prongs, and 1 singularity with 13 prongs.Ê The pseudo-Anosov acts as a five cycle on the 3-prong singularities.Ê ÊSo there are 2-singular orbits in the pseudo-Anosov flow on the 3-manifold.Ê ÊDeleting the 3-prong one gives something with homology $\mathbb Z \oplus \mathbb Z/5$, which is thus not a knot in $S^3$.Ê ÊDeleting the 13-prong one gives the census manifold $m011$, whose hyperbolic volume is 2.781833." Ê Aaber and Dunfield do not know the transition matrix (i.e. the digraph) for the pA flow, and indeed we just showed that it is only determined up to a choice of 5 different digraphs.  They also do not know the map that induces it explicitly, an interesting problem that we have not considered in this paper. This example should be interpreted together with the results in \cite{H}, where it is shown that for fixed $(d,a)$, particular transition matrices with LT polynomial $\xi_{d,a}$, are topologically induced, however Hironaka's methods would produce a pA map on a surface of genus 6 or 7, not 5.  The full picture is clearly not at hand. }

\begin{example} \label{Ex:palindromic but not pA} {\rm There are palindromic polynomials which cannot be induced by a pA mapping class acting on a surface of genus g.   For simple examples, let $T$ be a transition matrix of dimension $m$ for a pA map, and let 
$\xi_T(x) = x^m + \sum_{i=1}^m b_ix^{m-i}$.  Assume that $\xi_T(x)$ is palindromic.  Then $b_1 = $- trace$(T) \leq 0$ because $T$ is non-negative.   Therefore, for example, $x^2 +x +1$ cannot be the characteristic polynomial of $T$.  Another reason for the same phenomenon is that if $b_1>0$, then there would be a linear subdigraph of the associated digraph $D$ with 1 vertex, through which there passes an even number of disjoint cycles, which is impossible.   
 
A more interesting example is given by the infinite sequence of palindromic polynomials:
$$ p_n(x) = x^{24} +n(x^{20} - x^{19}) -x^{13} -x^{12} -x^{11} +n(-x^5 + x^4) +1,  \ n\geq 1$$
These polynomials are obtained from the Lanneau-Thiffeault polynomial $\xi_{12,11}$ by adding to it the polynomial
$$q_n(x)  = n(x^{20} - x^{19} -x^5 + x^4) =  n(x^{15} -1)x^4(x-1), $$
and since $q_n(1) =0$, also $q_n(x) > 0$ for all $x>1$, it follows that the largest real root of $p_n(x)$ is always bigger than 1, but approaches arbitrarily close to 1 as $n\to\infty$.  This cannot happen for infinitely many values of $n$, by \cite{AY} and \cite{Iv}.}
\end{example}

\begin{example} \label{Ex:the case n=1} 
{\rm We give an example of a digraph with antipalindromic characteristic polynomial, that illustrates the subtle points about the Coefficient Theorem for Digraphs. 
Figure 5 gives the digraph  for the monodromy map for the hyperbolic fibered knot $8_9$. 
\begin{figure}[htpb!]  \label{F:digraph4}
\centerline{\includegraphics[scale=.35] {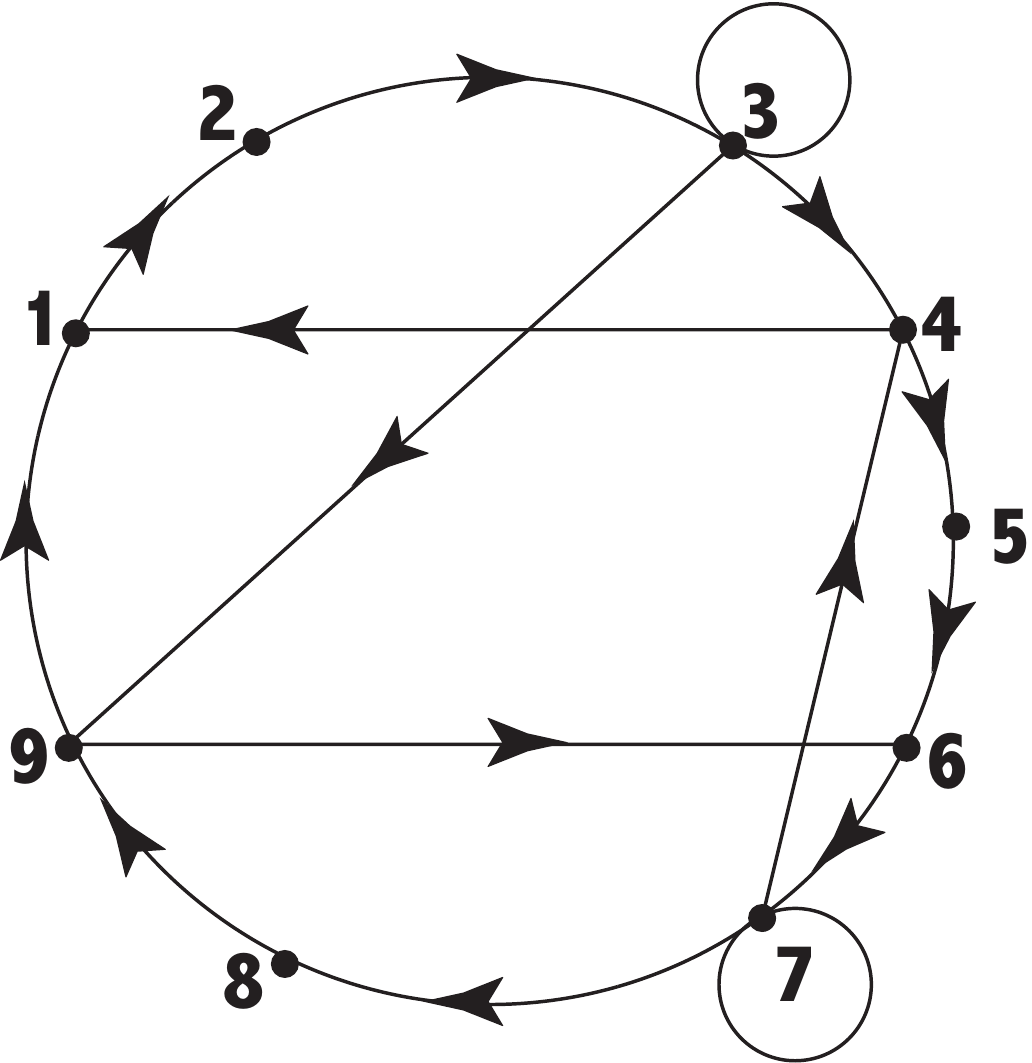}}
\caption{A digraph $D(1,4)$ that has an unexpected cycle} 
\end{figure}
A single  cycle  includes the 9 vertices and 9 connecting edges.  There are 6 additional edges, 2 of them loops based at vertices 3 and 7.   The strange feature of this example is the unexpected 7-cycle $L_7$ with vertices $1,2,3,9,6,7,4$.  While we did not find examples like this when we restricted to 4 or 5 additional edges,  this kind of surprise makes it difficult to use the Coefficient Theorem as $c-n$ increases. The reader should be able to check, easily, that the characteristic polynomial is 
$x^9 -2x^8 +x^7 -4x^5 +4x^4 -x^2 +2x -1$.}
\end{example} 

\begin{example}\label{Ex:bounds on c and d}
{\rm   We note that $\lambda_{d,a}$ decreases with increasing $d$, and increases with increasing $a$.  As an example, we consider the case $g=11$, and the difficulties presented by the polynomials in Equation (\ref{E:c=4 examples}).  We conjecture that they can be eliminated as candidates for $\lambda_{min,g}$ because,  with the restrictions on $\lambda_{min,g}$ that are given in  Equation (\ref{E:H}),  their spectral radius is too large.  But that is not true unless we can prove that if $\lambda$ is small, then the dimension of $T$ can't wander too far from $2g$.  

In our example the upper bound for $\lambda_{min,g}$ is the largest real root of $\xi_{12,11}(x) = x^{22}-x^{12}-x^{11}-x^{10}+1$, and Mathematica computes this to be $1.10918$. 
Now, we know that $2g\leq 2d \leq 6g-6$, or $11\leq d \leq 30$ in our example, and since we learned experimentally that $\lambda_{d,(a_1,a_2,a_3,a_4)}$ is smallest when $d$ is as large as possible and when the $a_i's$ are as close as possible to $d/2$, we were led to compute
the largest real root of $\xi_{30,(15,15,15,15)}= x^{60}-4x^{45}+5x^{30}-4x^{15}+1$, that is $\lambda_{30,(15,15,15,15)} = 1.06626.$  Since $1.06626 < 1.10918$, this tells us that we cannot omit consideration of the polynomials in Equation (\ref{E:c=4 examples}).  For comparison, we note that the largest real root of the LT polynomial $\xi_{30,29}(x)=x^{60}-x^{31}-x^{30}-x^{29}+1$ is $\lambda_{30,29} = 1.03262< 1.06626 = \lambda_{30,(15,15,15,15)}$, but of course it's very unlikely that this LT polynomial can be realized by a pA map on a surface of genus 11.  This implies that an additional measure of complexity that has been lurking in the background is $m-2g$, i.e. the difference between the dimension of an invariant  train track for a pA map and twice the genus of the surface on which it acts.  When $\lambda\leq \lambda_{d,a}$, one expects that $m-2g$ cannot be anything like $6g-8$, the only known upper bound.  In fact, in the all examples we know as candidates for $\lambda_{min,g}$  the quantity $m-2g$ is never bigger than 2. }
\end{example}

\begin{example} \label{Ex:Ham-Song Lemma} 
{\rm It was proved by Ham and Song in \cite{HS} that:
\begin{equation} \label{E:HS Lemma}
c \leq (\lambda_T)^m -1.
\end{equation}
The Ham-Song Lemma sounds like exactly what we need to know to prove that a given candidate for $\lambda_{min,g}$ is in fact $\lambda_{min,g}$, however the fly in the ointment is the exponent $m$ in Equation (\ref{E:HS Lemma}), the dimension of the transition matrix.   This brings us back to the conjecture in Example~\ref{Ex:bounds on c and d}.

The underlying fact that's used in the proof of the Ham-Song Lemma in \cite{HS} is that the Perron-Frobenius eigenvalue of an irreducible non-negative integer matrix is bounded below by the minimum and row sum in $T$, together with some elementary facts about powers of non-negative matrices.  That is not enough to give us the bound on $c$ that we would like to have.  We remark that  for genus 5 the dimension of $T$, that is the exponent in  (\ref{E:HS Lemma}), satisfies $10\leq m \leq 24$.  Since $\lambda_T > 1$, this exponent is crucial in determining the upper bound in (\ref{E:HS Lemma}).  The smallest known dilatation for genus 5  comes from an example in \cite{AD} whose transition matrix has dimension 14, not 24.  In all other examples the best known upper bounds have dimension $2g+2$.  If we could prove that for fixed $g$ and for pA maps with dilatation bounded above by $\Lambda_g$, the dimension $m$ of $T$ is at most $2g+4$, then the Ham-Song Lemma would tell us enough so that, with the help of a computer search, we could enumerate all possible cases for $\lambda_{d,a}$.   Nevertheless,  Equation~(\ref{E:HS Lemma}) suggests that low dilatation means low complexity. }
\end{example}


\noindent {\sc Department of Mathematics, \\
Barnard College and Columbia University,\\
2990 Broadway, New York, NY 10027}
\\ jb@math.columbia.edu


\begin{thebibliography}{{\bf BM}9}
\markright{}
\baselineskip18pt

\bibitem{AD}  J.W. Aaber and N. Dunfield, {\it Closed surface bundles of least volume}, preprint, arXiv:1002.3423.

\bibitem{AY} P.Arnoux and J.P. Yoccoz, {\it Construction de diff\'eomorphismes pseudo-Anosov}, C.R. Acad. Sci. Paris S\'er. I. Math. {\bf 292} (1981), no. 1, 75-78.

\bibitem{BH} M. Bestvina and M. Handel, {\it Train tracks for surface homeomorphisms}, Topology {\bf 34}, No. 1 (1995), 109-140.


\bibitem{BBK} J. Birman, P. Brinkmann and K. Kawamuro, {\it A polynomial invariant of pseudo-Anosov maps}, preprint,  arXiv:1001.5094.


\bibitem{CDS} D. Cvetkovic,  M. Doob and H. Sachs, {\it SPECTRA OF GRAPHS}, Academic Press 1980.  ISBN 0-12-195 150-2, LCCCN 79-50490.

\bibitem{FLM} B. Farb, C. Leininger and D. Margalit, {\it Small dilatation pseudo-Anosovs and 3-manifolds}, preprint,
arXiv:0905.0219

\bibitem{FLP} A. Fathi, F. Laudenbach and V. Poenaru, {\it TRAVAUX DE THURSTON SUR LES SURFACES}, Asterisque {\bf 66-67} (!979), Soci\'et\'e Math\'ematique de France, Paris.

\bibitem{HS} J.Y. Ham and W.T.Song, {\it The minimum dilatation of pseudo-Anosov 5-braids},  Journal of Experimental Mathematics {\bf 16}, No. 2 (2007), 167-179.   

\bibitem{H} E. Hironaka, {\it Small dilatation mapping classes coming from the simplest hyberbolic braid}, 
Algebraic and Geometric Topology, {\bf 10}, No. 4 (2010), 2041-2060.

\bibitem{Iv}  N. Ivanov 1988, {\it Stretching factors of pseudo-Anosov homeomorphisms},  Journal of Soviet Mathematics {\bf 52} (1990), 2819-2822.  English language translation of Zap. Nauchu. Sem. Leningrad Otdel. Mat. Inst. Steklov. (LOMI) {\bf 167} (1988), 111-116.


\bibitem{LT} E. Lanneau and J.L. Thiffeault, {\it On the minimum dilatation of pseudo-Anosov homeomorphisms on surfaces of small genus}, preprint, arXiv:0905.1302.

\bibitem{McM} C. McMullen, {\it Polynomial invariants for fibered 3-manifolds and Teichmuller geodesic for foliations}, Annales Scientifique de l'\'Ecole Normal Sup\'erieure, Quatri\`eme S\'erie {\bf 33} (2000), 519-560.


\bibitem{R} I. Rivin, {\it  Walks on groups, counting reducible matrices, polynomials, and surface and free group automorphisms}, Duke Math. Journal {\bf 142}, No. 2 (2008), 353-379.

\bibitem{Se}  Seneta, {\it NON-NEGATIVE MATRICES and MARKOV CHAINS}, Springer Science and Business Series, 1973, 1981 and 2006.

\bibitem{Th}
W.P. Thurston, {\em On the geometry and dynamics of diffeomorphisms of surfaces}, Bull. Amer. Math. Soc., {\bf 19}, No.2 (1988), 417--431.


\end{thebibliography}
\end{document}